# FRACTIONAL ISOPERIMETRIC INEQUALITIES AND SUBGROUP DISTORTION

Martin R. Bridson

ABSTRACT. It is shown that there exist infinitely many non-integers $r > 2$ such that the Dehn function of some finitely presented group is $\simeq n^r$. For each positive rational number $s$ we construct pairs of finitely presented groups $H \subset G$ such that the distortion of $H$ in $G$ is $\simeq n^s$. And for each $s \geq 1$ we also construct finitely presented groups whose isodiametric function is $\simeq n^s$.

*For John Stallings, on his 60 th birthday*

**Introduction.**

Isoperimetric inequalities measure the complexity of the word problem in finitely presented groups by giving a bound on the number of conjugates of relators that one must apply in order to show that a word $w$ in the given generators represents the identity. Such bounds are given in terms of the length of $w$, and the function describing the optimal bound is known as the *Dehn function* of the group.

About four years ago, as the result of efforts by a number of authors [BMS], [BP], [Gr2], it was established that for every positive integer $d$ one can construct finitely presented groups whose Dehn function is polynomial of degree $d$. The question of whether or not there exist groups for which the optimal isoperimetric inequality is of fractional degree has attracted a good deal of interest (e.g., [N, prob. 32], [SS]). According to a theorem of M. Gromov, there are no such fractional exponents less than 2, because if a group satisfies a sub-quadratic isoperimetric inequality then it actually satisfies a linear isoperimetric inequality (see [Gr1], [Ol], [Pa]).

The main purpose of this article is to prove the following:

**Theorem A.** *There exist infinitely many non-integers $r > 2$ such that the Dehn function of some finitely presented group is $\simeq n^r$.*

Our method for constructing groups whose Dehn functions are of fractional degree is geometric in nature. Of the examples which we shall describe, the simplest are obtained by taking three torus bundles over the circle (each of a different dimension) and amalgamating their fundamental groups over central cyclic subgroups.

Given a smooth closed curve in the universal cover of a compact Riemannian manifold $M$, the isoperimetric function of $M$ bounds the area of a minimal spanning disc which fills the loop; the bound is given as a function of the length of the loop. One can show that the isoperimetric function of $M$ is $\simeq$ equivalent to the Dehn

This work was supported in part by NSF grant DMS-9203500

Typeset by $\mathcal{A}\mathcal{M}\mathcal{S}$-TEX

function of its fundamental group. Thus Theorem A can be interpreted as a result about isoperimetric functions of Riemannian manifolds.

Lurking in the background of our proof of Theorem A is an observation concerning subgroup distortion. Following Gromov [Gr2], we define the distortion of a pair of finitely generated groups $H \subseteq G$ to be the $\simeq$ equivalence class of the function $\delta : \mathbb{N} \to \mathbb{N}$, where $n\, \delta(n)$ is the radius of the set of vertices in the Cayley graph of $H$ that are a distance at most $n$ from the identity in $G$.

**Theorem B.** *For every positive rational number $r$ there exist pairs of finitely presented groups $H \subseteq G$ with distortion $\simeq n^r$.*

In certain contexts, instead of bounding the number of conjugates of relators that one must apply in order to show that a word represents the identity in a group $G$, it is more appropriate to bound the length of the words conjugating the relators. The optimal such bound is given by the *isodiametric function* of $G$. If $G$ is the fundamental group of a compact Riemannian manifold $M$, then this function gives a bound on the diameter of spanning discs in $M$.

**Thereom C.** *For every positive rational number $r \geq 1$ there exist finitely presented groups whose isodiametric function is $\simeq n^r$.*

In [B4] we shall explain how Thereoms A and C can be used to answer a question of M. Gromov concerning the connection between the Dehn functions of finitely presented groups and the topology of their asymptotic cones.

The results contained in this paper were described in seminars in Princeton and New York at the end of 1993 and more fully in my lectures at the conferences on geometry and group theory in Champousin, Durham and Lyon in the spring and summer of 1994. I would like to thank the organisers of all of these excellent conferences. I would also like to apologise for my tardiness in preparing the final draft of this paper. In the meantime, further significant progress on the fundamental question of which functions arise as Dehn functions of finite presentations has been made by J-C. Birger, E. Rips and M. Sapir [BRS]. Although the present extent of their results is unclear, it appears that by encoding descriptions of certain turing machines into finite presentations they can in particular construct finitely presented groups for which the Dehn function has the form $n^r$ with $r$ irrational.

**Section 1: Definitions.**

We begin by recalling the definition of the Dehn function of a finite presentation. We fix a finite presentation $\mathcal{P} = \langle \mathcal{A} \mid \mathcal{R} \rangle$ for the group $\Gamma$. We assume that $\mathcal{R}$ is symmetrized [LS]. A reduced word $w \in F(\mathcal{A})$ is *null-homotopic* (i.e., represents the identity element in $\Gamma$) if and only if there is an equality of the form $w = \prod_i^N g_i^{-1} r_i g_i$ in the free group $F(\mathcal{A})$, where $r_i \in \mathcal{R}$. Isoperimetric inequalities give bounds on the number $N$ of factors in a minimal such expression; this integer is called $Area(w)$. One seeks to bound $Area(w)$ in terms of $|w|$, the length of $w$. The function giving the optimal such bound is called the Dehn function of the presentation.

For the purposes of this article, it is best to cast the definition of the Dehn function in the following geometric context.

According to van Kampen's lemma (see [LS]) the above factorisation of the null homotopic word $w$ can be portrayed by a finite, oriented, planar, combinatorial 2-complex with basepoint. This complex is called a van Kampen diagram for $w$. The oriented 1-cells of this complex are labelled by elements of $\mathcal{A}$ and their inverses, the boundary label on each face of the 2-complex is an element of $\mathcal{R}$, and the boundary cycle of the complex (read with positive orientation from the basepoint) is the word $w$. The number of 2-cells in the complex is equal to $N$, the number of factors in the given equality for $w$. Conversely, any van Kampen diagram gives rise to an equality in $F(\mathcal{A})$ showing that the boundary cycle of the diagram represents the identity in $\Gamma$. (See [BG] for a more precise treatment.) A *minimal area* van Kampen diagram is one that has the least number of 2-cells among all van Kampen diagrams which share its boundary label. Notice that a connected and simply connected subdiagram of a minimal area diagram again has minimal area.

*1.1 Definition.* If $\mathcal{P}$ is a finite presentation for the group $\Gamma$, and if $w$ is a word in the generators and their inverses that represents the identity in $\Gamma$, then $Area(w)$ is defined to be the number of 2-cells in a minimal area van Kampen diagram for $w$. The *Dehn function* $f : \mathbb{N} \to \mathbb{N}$ *of* $\mathcal{P}$ is defined by:

$$f(n) = \max_{|w| \leq n} Area(w).$$

If $\mathcal{Q}$ is another finite presentation for $\Gamma$ with Dehn function $g$, then there exist positive constants $A, B, C, D, E$ so that for all $n \geq 0$ [A]. All of the Dehn functions which we shall consider in this article grow at least linearly (and we shall implicitly assume this henceforth). In this setting one can dispense with the term $(Dn + E)$ in $(*)$. This motivates a notion of equivalence of functions, which we now recall.

(*) $$g(n) \leq Af(Bn + C) + Dn + E,$$

In general, given $f, g : \mathbb{N} \to \mathbb{N}$, one writes $g \preceq f$ if there exist positive constants $A, B, C$ so that $(*)$ holds with $D = E = 0$. And one writes $f \simeq g$ if, in addition, $f \preceq g$. It is easily verified that this is an equivalence relation. Henceforth we shall speak of "the Dehn function of the group $\Gamma$" with the understanding that this is only defined up to $\simeq$ equivalence.

Throughout this article we shall be concerned with presentations of the following type. Let $\mathcal{A}_1, \ldots, \mathcal{A}_n$ be a disjoint collection of finite sets, and consider:

(1.2) $$\langle \mathcal{A}_1 \cup \cdots \cup \mathcal{A}_n, z_1, \ldots, z_{n-1} \mid \mathcal{R}_1, \ldots, \mathcal{R}_n \rangle$$

where the relators $\mathcal{R}_i$ involve only the letters $\mathcal{A}_i \cup \{z_i, z_{i-1}\}$, if $i \in \{2, \ldots, n-1\}$, the relators $\mathcal{R}_1$ involve only the letters $\mathcal{A}_1 \cup \{z_1\}$, and the relators $\mathcal{R}_n$ involve only the letters $\mathcal{A}_n \cup \{z_{n-1}\}$. Assume that in the group defined by this presentation the elements $z_i$ all have infinite order. Notice that one encounters presentations such as this when considering groups of the following form: $G_1 *_{C_1} \cdots *_{C_{n-1}} G_n$, where the amalgamated subgroups $C_i$ are all infinite cyclic.

Given a van Kampen diagram $\Delta$ over such a presentation, we consider the equivalence relation $\sim$ on 2-cells in $\Delta$ generated by setting $e \sim e'$ if the boundary cycle of each of $e$ and $e'$ is labelled with a relator from the same $\mathcal{R}_i$ and if, in addition, their boundaries intersect in at least one edge.

*1.3 Definition.* A *monochromatic region* of $\Delta$ is the union of an $\sim$ equivalence class of closed 2-cells. If the boundary labels of the 2-cells in this equivalence class belong to $\mathcal{R}_i$, then the monochromatic region is said to be of type $\mathcal{R}_i$.

For example, in the case $n = 1$, the monochromatic regions are precisely the discs obtained as the closures of the connected components of $\Delta$ minus its boundary cycle $\partial \Delta$.

**1.4 Lemma.** *If $\Delta$ is minimal area diagram then all of its monochromatic regions are homeomorphic to discs, and every monochromatic region of type $\mathcal{R}_1$ or $\mathcal{R}_n$ meets the boundary of $\Delta$ in at least one edge.*

*Proof.* We proceed by induction on $n$. As the preceding remark indicates, the case $n = 1$ is clear.

If $\Delta$ did not contain a monochromatic region of type $\mathcal{R}_1$ then we would be done by application of our inductive hypothesis to $\Delta$ viewed as a van Kampen diagram over the sub-presentation obtained by deleting $\mathcal{A}_1$ and $\mathcal{R}_1$ and incorporating $z_1$ into $\mathcal{A}_2$.

We proceed under the assumption that there is a monochromatic region of type $\mathcal{R}_1$.

By definition, the frontier of a region is the intersection of the region with the closure of its complement in the plane. We claim that the frontier of each monochromatic region of type $\mathcal{R}_1$ is homeomorphic to a circle (implying that the region is a disc) and that, moreover, it also contains an edge in the boundary of $\Delta$. Indeed, if either of these two conditions were to fail then the frontier of our region would contain a simple closed loop all of whose 1-cells were in the interior of $\Delta$. But if a 2-cell of type $\mathcal{R}_1$ meets a 2-cell of another type along an edge, then the second 2-cell must be of type $\mathcal{R}_2$ and the edge must be labelled $z_1$. Hence this simple closed loop in the interior of $\Delta$ would be labelled by a word involving only $z_1$ and $z_1^{-1}$. But $z_1$ is assumed to have infinite order in the group defined by the above presentation, so this labelling word freely would reduce to the trivial word, i.e., have area 0. But then, by replacing the (simply connected) subdiagram of $\Delta$ bounded by this loop with a diagram of area 0, we would have a contradiction to the hypothesis that $\Delta$ is of minimal area.

Having established the asserted properties for regions of type $\mathcal{R}_1$ (and hence, by symmetry, type $\mathcal{R}_n$) it only remains to observe that, given a disc subdiagram of a van Kampen diagram, if the subdiagram contains an edge of the boundary of the ambient diagram, then the closure of its complement is a disjoint union of van Kampen diagrams. Thus we may remove all monochromatic regions of type $\mathcal{R}_1$ from $\Delta$, and applying our inductive hypothesis to the resulting collection of van Kampen diagrams (considered as diagrams over the sub-presentation described in the first paragraph of the proof) we deduce that all monochromatic regions of $\Delta$ are homeomorphic to discs. $\square$

*1.5 Remark.* The above argument remains valid if the cyclic groups $\langle z_i \rangle$ are replaced by free groups of finite rank.

**1.6 Definition.** Let $G$ be a group with finite generating set $\mathcal{A}$. A word $w$ in the free group on $\mathcal{A}$ is called *injective* if no non-empty subword of $w$ is null-homotopic (i.e., represents $1 \in G$).

A null-homotopic word $v$ in the free group on $\mathcal{A}$ is called *irreducible* if no proper subword of $v$ is null-homotopic. (This implies that every van Kampen diagram for $v$ is a topological disc.)

**Section 2: The groups $G_c$ and $\Gamma(a,b,c)$.**

The Dehn functions for many groups of the form $\mathbb{Z}^m \rtimes \mathbb{Z}$ were calculated in [BP], and this classification was later extended to all abelian-by-free groups (see [BG], [B3]). We shall be concerned with a particular class of examples from [BP], namely the groups $G_c = \mathbb{Z}^c \rtimes_{\phi_c} \mathbb{Z}$, where $\phi_c \in Gl_c(\mathbb{Z})$ is the unipotent matrix with ones on the diagonal and super-diagonal and zeros elsewhere. $G_c$ has presentation:

$$\langle x_1, \ldots, x_c, t \mid [x_i, x_j] = 1 \text{ for all } i,j,\ [x_c, t] = 1,\ [x_i, t] = x_{i+1} \text{ if } i < c \rangle.$$

Notice that the centre of $G_c$ is the infinite cyclic subgroup generated by $x_c$. To emphasize this fact, we shall change notation and write $z_c$ in place of $x_c$.

In what follows we shall use $d$ to denote the word metric on $G_c$ corresponding to the generators given in the above presentation. (Recall that the word metric associated to a finite system of generators for a group $G$ is the left-invariant metric in which the distance of each element $g \in G$ from the identity is the length of the shortest word in the generators and their inverses that represents $g$.)

We shall need the following facts concerning $G_c$.

**2.1 Proposition.**
   (1) *The Dehn function of $G_c$ is $\simeq n^{c+1}$.*
   (2) *There exists a constant $\alpha_c$ such that $Area(wz^n) \leq \alpha_c |w|^{c+1}$ for all $n \in \mathbb{Z}$ and all words $w$ such that $wz^n = 1$ in $G_c$.*
   (3) *There exist constants $K_c > 1$ and $k_c > 0$ such that, for all $n > 0$,*

$$k_c n^{1/c} \leq d(1, z_c^n) \leq K_c n^{1/c}.$$

*Proof.* Assertion (1) is proved in [BP]. The way in which upper bounds are obtained on Dehn functions in that article is by using the combings of $G_c$ constructed in [B2]; one then appeals to results of [B1] to see that the area of a word (with respect to a certain finite presentation) is bounded by the product of the length of the word and the length of the longest combing line to any of its vertices. (2) is an easy consequence of the fact that $z^n$ is a combing line for the (asynchronously bounded) combing of $G_c$ which was constructed in [B2].

(3) is an observation of Gromov [Gr2, $5.A_2'$] (see also [Pi]). The key point is that the Lie algebra of the Lie group $\mathbb{R}^c \rtimes_{\phi_c} \mathbb{R}$ is graded, and hence the Lie group admits a 1-parameter family of automorphisms $H_t$, with $t \in [0,1]$, such that for all $x, y \in \mathbb{R} \rtimes_{\phi_c} \mathbb{R}$, one has $d_{CC}(H_t(x), H_t(y)) = t\, d_{CC}(x,y)$, where $d_{CC}$ is the Carnot-Carathéodory metric on $\mathbb{R}^c \rtimes_{\phi_c} \mathbb{R}$ associated to the field of tangent 2–planes orthogonal to the kernel of the derivative of the abelianization map $\mathbb{R}^c \rtimes_{\phi_c} \mathbb{R} \to \mathbb{R}^2$. For each element $g = \exp(u)$ of the centre, $H_t(g^c) = \exp(tc.u)$. Hence, setting $\tau = 1/c$, one gets: $d_{CC}(1, g^c) = \tau\, d_{CC}(1, \exp(\tau c.u)) = (1/c)\, d_{CC}(1, g)$. The desired inequality for $G_c$ is then an immediate consequence of the fact that $G_c$ is a cocompact lattice in $\mathbb{R} \rtimes_{\phi_c} \mathbb{R}$, and hence is quasi-isometric to $\mathbb{R} \rtimes_{\phi_c} \mathbb{R}$ equipped with the (left-invariant) Carnot-Carathéodory metric. □

We shall consider the groups $\Gamma(a, b, c)$ obtained as follows: first we amalgamate $G_a$ with $G_b \times \mathbb{Z}$ by identifying the centre of $G_a$ with that of $G_b$; then we form the amalgamated free product of the resulting group with $G_c$, identifying the centre of the latter with the right-hand factor of $G_b \times \mathbb{Z}$. In symbols:

$$\Gamma(a, b, c) = G_a *_{z_a = z_b} (G_b \times \langle \zeta \rangle) *_{\zeta = z_c} G_c.$$

Combining the natural presentations for the $G_l$ we obtain:

**2.2 Presentation of $\Gamma(a, b, c)$.**
*Generators:*

$$x_1, \ldots, x_{a-1}, z(= z_a = z_b), t_a, y_1, \ldots, y_{b-1}, t_b, u_1, \ldots, u_{c-1}, \zeta(= z_c), t_c.$$

*Defining relations:*

$$[z, x_i] = [x_i, x_j] = 1 \ \forall i, j; \ [z, t_a] = 1; \ [x_i, t_a] = x_{i+1} \ if \ i < a-1 \ and \ [x_{a-1}, t_a] = z;$$

*and*

$$[z, y_i] = [\zeta, y_i] = [y_i, y_j] = 1 \ \forall i, j; \ [z, t_b] = [\zeta, t_b] = [z, \zeta] = 1;$$

$$[y_i, t_b] = y_{i+1} \ if \ i < b-1, \ and \ [y_{b-1}, t_b] = z;$$

*and*

$$[\zeta, u_i] = [u_i, u_j] = 1 \ \forall i, j; \ [\zeta, t_c] = 1; \ [u_i, t] = u_{i+1} \ if \ i < c-1 \ and \ [x_{c-1}, t-c] = \zeta.$$

Let $\mathcal{A}_a$ denote the generators $x_i$ together with $t_a$, let $\mathcal{A}_b$ denote the generators $y_i$ together with $t_b$, and let $\mathcal{A}_c$ denote the generators $u_i$ together with $t_c$. Let $\mathcal{R}_a$ denote the first set of relators, $\mathcal{R}_b$ the second set and $\mathcal{R}_c$ the third set.

Notice that the above presentation $\langle \mathcal{A}_a, \mathcal{A}_b, \mathcal{A}_c, z, \zeta \mid \mathcal{R}_a, \mathcal{R}_b, \mathcal{R}_c \rangle$ is precisely the sort of presentation to which (1.3) and (1.4) apply. By exploiting this fact we shall prove:

**2.3 Theorem.** *For all integers $1 \leq b \leq a < c$, the Dehn function of $\Gamma(a, b, c)$ is $\simeq n^{c + \frac{a}{b}}$.*

**Section 3: Lower bounds and subgroup distortion.**

The purpose of this section is to prove the following:

**3.1 Proposition.** *For any positive integers $a, b$ and $c$, the Dehn function $f(n)$ for the above presentation of $\Gamma(a, b, c)$ satisfies*

$$n^{c+\frac{a}{b}} \preceq f(n).$$

The clarity of our exposition will be enhanced if we pursue some general considerations first.

**3.2 Definition of distortion.** Let $H \subset G$ be a pair of finitely generated groups, and let $d_G$ and $d_H$ be the word metrics associated to a fixed choice of generators for each. Following Gromov [Gr2], we define the *distortion* of $H$ in $G$ to be the function

$$\delta(n) = \frac{1}{n} \max\{d_H(1, h) \mid h \in H \text{ with } d_G(1, h) \le n\}.$$

One checks easily that, up to $\simeq$ equivalence, this function is independent of the choice of word metrics $d_G$ and $d_H$.

**3.3 Example.** According to 2.1(3), the distortion of the centre $Z(G_c)$ in $G_c$ is $\simeq n^{\frac{1}{c}-1}$ for all $c \in \mathbb{Z}_+$. And by comparing $d_{G_a}(1, z^n)$ with $d_{G_b}(1, z^n)$, we see that if $a > b$ then the distortion of $G_b$ in $G_a *_{\langle z \rangle} G_b$, the group formed by amalgamating $G_a$ and $G_b$ along their centres, is $\succeq n^{\frac{a}{b}-1}$.

If $G$ is a group with finite presentation $\langle \mathcal{A} \mid \mathcal{R} \rangle$ and $H$ is the subgroup generated by a subset $\mathcal{B} \subset \mathcal{A}$ then the HNN extension of $G$ obtained by adding a stable letter that commutes with $H$ has finite presenatation

$$\mathcal{P} = \langle \mathcal{A}, \tau \mid \mathcal{R}, [\tau, b] = 1 \; \forall b \in \mathcal{B} \rangle.$$

We shall now explain how the distortion of $H$ in $G$ provides a lower bound on the Dehn function of this HNN extension.

**$\tau$-corridors.** Consider a reduced van Kampen diagram $\Delta$ over $\mathcal{P}$, and consider an oriented edge of $\partial \Delta$ that is labelled $\tau$. If this edge lies in the boundary of some 2-cell then the relation labelling this 2-cell must be of the form $\tau b \tau^{-1} b^{-1}$. In particular, the 2-cell has a single oriented edge labelled $\tau$ besides the one that we started with. If this second $\tau$-edge is in the interior of $\Delta$, then a second 2-cell with a boundary label of the form $\tau b' \tau^{-1} b'^{-1}$ must contain this edge in its boundary. (And because the diagram is reduced, $b' \ne b^{-1}$.) By iterating this argument one obtains a chain of 2-cells crossing $\Delta$, beginning at the $\tau$-edge that we first considered and ending at some other $\tau$-edge in $\partial \Delta$. Such a chain of 2-cells is called a *$\tau$-corridor*.

A rigorous treatment of $\tau$-corridors (in greater generality) can be found in [BG]. In the present setting we need only observe that $\tau$-corridors exist, that no 2-cell of $\Delta$ lies in more than one $\tau$-corridor, that the 'sides' of a $\tau$-corridor are labelled by words in the generators $\mathcal{B}$, and that the area of a $\tau$-corridor is equal to the length of each of its sides. (More precisely, each $\tau$-corridor is a subdiagram of $\Delta$ whose boundary cycle is labelled $\tau v \tau^{-1} v^{-1}$, where $v$ is a word in the generators $\mathcal{B}$, and the number of 2-cells in the corridor is $|v|$.)

**3.4 Proposition.** *If $G$ is a finitely presented group and $H$ is a finitely generated subgroup with distortion $\delta(n)$, then the Dehn function $f(n)$ of the HNN extension $\{G, \tau : [\tau, h] = 1 \ \forall h \in H\}$ satisfies:*

$$n\,\delta(n) \preceq f(n).$$

*More generally, if $G = \langle \mathcal{A} \mid \mathcal{R} \rangle$ and $H$ is the subgroup generated by $\mathcal{B} \subset \mathcal{A}$, and if $v$ is an injective word in the generators $\mathcal{A}$ that represents $h \in H$, then for all $m \in \mathbb{N}$, with respect to the presentation $\langle \mathcal{A}, \tau \mid \mathcal{R}, [\tau, b] = 1 \ \forall b \in \mathcal{B}\rangle$,*

$$m\,d_\mathcal{B}(1, h) \leq \mathrm{Area}(v\tau^m v^{-1} \tau^{-m}).$$

*Proof.* Because $v$ is injective, $W = v\tau^m v^{-1}\tau^{-m}$ is irreducible. In particular, any van Kampen diagram for $W$ is a disc, so each edge in its boundary lies in the boundary of some closed 2-cell. We fix a reduced diagram $\Delta$ for $W$.

There is a $\tau$-corridor beginning at each edge in the segment of $\partial\Delta$ labelled $\tau^m$. The $\tau$-corridor beginning at the $i$-th edge of the segment of $\partial\Delta$ labelled $\tau^m$ ends at the $(m-i)$-th edge of the segment of $\partial\Delta$ labelled $\tau^{-m}$, and the side of each such corridor is labelled by a word in the generators $\mathcal{B}$ and their inverses that represents the group element $h$. Thus each $\tau$-corridor contains at least $d_\mathcal{B}(1, h)$ 2-cells. $\square$

**Proof of Proposition 3.1.** We work with respect to the presentation of $\Gamma(a, b, c)$ defined in the previous section. Let $U_n$ denote the shortest word in the generators $\mathcal{A}_c$ and their inverses that represents $\zeta^{n^c}$, and let $V_n$ denote the shortest word in the generators $\mathcal{A}_a$ and their inverses that represents $z^{n^a}$. Because $z = [x_{a-1}, t_a]$, the length of $V_n$ is less than $4\,d_{G_a}(1, z^{n^a})$, which by 2.1(3) is no more than $4K_a n$. Similarly $|U_n| \leq 4K_c n$. Thus the family of words $W_n = U_n V_n U_n^{-1} V_n^{-1}$ have lengths $|W_n| \simeq n$. Because $z$ and $\zeta$ commute, $W_n$ represents the identity in $\Gamma(a, b, c)$. We will show that $\mathrm{Area}(W_n) \succeq n^{c+\frac{a}{b}}$.

Let $\Delta$ be a minimal area van Kampen diagram for $W_n$. Because no proper subword of $W_n$ represents the identity in $\Gamma(a, b, c)$, $\Delta$ is a disc. Let $\sigma_1, \sigma_2, \sigma_3, \sigma_4$ denote the segments of $\partial\Delta$ labelled $U_n, V_n, U_n^{-1}, V_n^{-1}$ respectively. Let $p_1, p_2, p_3, p_4$ be the initial vertices of these arcs. We shall argue, using the monchromatic regions of Lemma 1.4, that there is an arc labelled $\zeta^{n^c}$ in the 1-skeleton of $\Delta$ that connects $p_1$ to $p_2$.

Of the two edges in $\partial\Delta$ incident at $p_1$, one is labelled by a letter from $\mathcal{A}_a$ and the other is labelled by a letter for $\mathcal{A}_c$. The former must lie in the boundary of a 2-cell of type $\mathcal{R}_a$ and the latter must lie in the boundary of a 2-cell of type $\mathcal{R}_c$. It follows that $p_1$ is the endpoint of an arc in the interior of $\Delta$ that bounds a monochromatic region of type $\mathcal{R}_c$. Let $q \in \partial\Delta$ be the other endpoint of this arc. Since the arc is labelled by a power of $\zeta$, the portion of $\partial\Delta$ connecting $p_1$ to $q$ must be labelled by a word that represents an element of $\langle \zeta \rangle$. Thus $q \in \sigma_1$. If $q \neq p_2$ then consider the edge of $\sigma_1$ that is between $q$ and $p_2$ and is incident at $q$. This edge lies in a monochromatic region of type $\mathcal{R}_c$, and emanating from $q$ there is an arc in the interior of $\Delta$ that bounds this monochromatic region. The other endpoint of this arc is a point $q' \in \sigma_1$ that is between $q$ and $p_2$ and distinct from $q$. Proceeding in this manner we can connect $p_1$ to $p_2$ by a concatenation of arcs all of whose edges are in the interior of $\Delta$ and are labelled by $\zeta$ or $\zeta^{-1}$. Because $\Delta$ has minimal area, its 1-skeleton does not contain any closed loops with such labelling, so this concatenation is actually a simple arc. In the same way we obtain a $\zeta$-labelled simple arc from $p_3$ to $p_4$.

The two paths which we have just constructed must be disjoint, because otherwise there would be a $\zeta$-labelled path from $p_1$ to $p_4$, contradicting the fact that the label on the arc of $\partial\Delta$ connecting $p_4$ to $p_1$ represents a non-trivial element of $\langle z \rangle$. Thus $\Delta$ contains a disc subdiagram $\Delta_0$ bounded by the two $\zeta$-labelled paths constructed above and the arcs of $\partial\Delta$ labelled $V_n$ and $V_n^{-1}$. This subdiagram cannot contain any 2-cells corresponding to relations of type $\mathcal{R}_c$, because this would contradict the fact (Lemma 1.4) that every monochromatic region of type $\mathcal{R}_c$ in $\Delta$ meets $\partial\Delta$ in at least one edge. Thus $\Delta_0$ can be viewed as a van Kampen diagram for $\zeta^{n^c} V_n \zeta^{-n^c} V_n^{-1}$ over the sub-presentation $\langle \mathcal{A}_a, \mathcal{A}_b, z, \zeta \mid \mathcal{R}_a, \mathcal{R}_b \rangle$.

This sub-presentation defines the subgroup $G_a *_z (G_b \times \langle \zeta \rangle)$ and the only relations involving $\zeta$ are $[z, \zeta] = 1$ and $[y_i, \zeta] = 1$ for all $y_i \in \mathcal{A}_b$. Thus we are in the situation of Proposition 3.4, with $\mathcal{A} = \mathcal{A}_a \cup \mathcal{A}_b \cup \{z\}$, with $\mathcal{B} = \mathcal{A}_b \cup \{z\}$ and with $\tau = \zeta$. By combining 2.1(3) and (3.4) we deduce:

$$k_b n^{c+\frac{a}{b}} \leq n^c \, d_{G_b}(1, z^{n^a}) \leq \text{Area}\,\Delta_0 \leq \text{Area}(W_n).$$

$\square$

## Section 4: Upper Bounds.

In this section we establish the upper bound necessary to complete the proof of Theorem 3.1. The most common method for obtaining upper bounds on Dehn functions is to construct a bounded (perhaps asynchronous) combing for the group and then estimate its length (in the sense of [B1]). This method is doomed to failure in our case, because the words $W_n$ considered in the preceding proof do not bound van Kampen diagrams whose diameter is bounded by a linear function of $|W_n|$, and this implies that the groups $\Gamma(a, b, c)$ are not asynchronously combable. Thus we are forced to undertake a direct analysis of van Kampen diagrams over the presentation of $\Gamma(a, b, c)$ given in Section 2.

We wish to establish an upper bound of the form $kn^r$ on the area of null-homotopic words. Because this function is superadditive (i.e. $k(s+t)^r \geq ks^r + kt^r$ for all $s, t > 0$), it suffices to consider null-homotopic words $w$ which are irreducible; for if $w = uv$ with $u$ and $v$ null-homotopic, then by forming the one-point join of minimal area diagrams for $u$ and $v$ we have $\text{Area}(w) \leq \text{Area}(u) + \text{Area}(v)$. As we noted in (1.6), the advantage of restricting our attention to irreducible words is that any van Kampen diagram for such a word is a topological disc.

> *We fix a word $w$ of length $n$ that is null-homotopic and irreducible.*
> *And we fix a minimal area van Kampen diagram $\Delta$ for $w$.*

We shall replace $\Delta$ by a cellulated disc whose combinatorial structure is simpler. In order to do so we focus on the monochromatic regions of type $\mathcal{R}_a$ and $\mathcal{R}_c$ in this diagram. The frontier of a monochromatic region of type $\mathcal{R}_a$ is topologically a circle (1.4). This circle contains a number of subarcs of $\partial\Delta$ (called $\partial$-arcs of type $\mathcal{R}_a$), and these are connected by arcs in the interior of $\Delta$ whose edges are labelled $z^{\pm 1}$; we call these connecting arcs *z-arcs*. The boundary of each monochromatic region of type $\mathcal{R}_c$ admits a similar description, but with $\zeta$-arcs in place of *z*-arcs. The combinatorial geometry of this system of labelled arcs encodes most (but not all) of the decomposition of $\Delta$ into monochromatic regions. (See figure 1.)

We now leave $\Delta$ behind and begin working with the following labelled cell complex structure $\Delta'$ on the unit disc in $\mathbb{R}^2$: the boundary circle is cellulated and labelled as the boundary of $\Delta$; we then focus on the vertices corresponding to the endpoints of $\partial$-arcs of type $\mathcal{R}_a$ and $\mathcal{R}_c$, and divide the interior of the disc into faces (2-cells) by connecting these vertices with non-intersecting, labelled arcs in the interior of the disc — two vertices on $\partial \Delta'$ are connected by an arc labelled $z^r$ (respectively, $\zeta^s$) if the corresponding points on $\partial \Delta$ were connected by a *z*-arc (respectively, a $\zeta$-arc) in $\Delta$ with the same label. (It is convenient to retain the terminology *z*-arc and $\zeta$-arc for the arcs just constructed in $\Delta'$.)

$\Delta'$ is a somewhat simplified version of $\Delta$. Each monochromatic region $D \subseteq \Delta$ of type $\mathcal{R}_a$ or $\mathcal{R}_c$, determines a face in $\Delta'$ whose boundary label is the same as that of $D$; in particular, it makes sense to speak of these faces in $\Delta'$ as being of type $\mathcal{R}_a$ or $\mathcal{R}_c$, and henceforth we shall do so. The remaining faces of $\Delta'$ are said to be of type $\mathcal{R}_b$, but we should immediately point out that there is not a 1-1 correspondence between monochromatic regions of type $\mathcal{R}_b$ in $\Delta$ and faces of type $\mathcal{R}_b$ in $\Delta'$: the word labelling the boundary of a face of $\Delta'$ that is of type $\mathcal{R}_b$ is a null-homotopic word in the generators of $G_b \times \langle \zeta \rangle$, and this word appears as the label on a circuit in the 1-skeleton of $\Delta$, but this word does not necessarily appear as the label on the boundary cycle of a single monochromatic region of type $\mathcal{R}_b$ in $\Delta$ (see figure 1).

Figure 1: Replacing $\Delta$ with $\Delta'$

In order to complete the proof we shall show that the sum of the areas of the words labelling the faces of $\Delta'$ is bounded by a suitable constant times $n^{c+(a/b)}$. To enhance the clarity of our exposition, we break the proof into a number of steps.

*Let the constants $\alpha_a, \alpha_b, \alpha_c, k_a, k_b, k_c, K_a, K_b, K_c$ be as in Proposition 2.1, and let*

$$K := \max\left\{K_a, K_b, K_c, \frac{1}{k_a}, \frac{1}{k_b}, \frac{1}{k_c}, \alpha_a, \alpha_b, \alpha_c\right\}.$$

*4.0 Remark.* In what follows we shall repeatedly use the hypothesis $a \geq b$ in order to obtain estimates of the form $x^{a/b} + y^{a/b} \leq (x+y)^{a/b}$.

**Step 1:** *Bounding the length of the internal arcs in $\Delta'$.*

**4.1 Lemma.** *If the $\zeta$-arcs in the boundary of a face in $\Delta'$ of type $\mathcal{R}_c$ are labelled $\zeta^{r_1}, \ldots, \zeta^{r_s}$, then $\sum_i |r_i| \leq Km^c$, where $m$ is the number of occurences of letters from $\mathcal{A}_c^{\pm 1} \cup \{\zeta^{\pm 1}\}$ in $w$ (the boundary label of $\Delta'$).*

*Proof.* Killing the generators of $\Gamma(a,b,c)$ other than $\mathcal{A}_c \cup \{\zeta\}$ gives a retraction of $\Gamma(a,b,c)$ onto $G_c$. So if an arc of $\partial \Delta'$ is labelled by a word $v$ that represents $\zeta^{r_i}$ in $\Gamma(a,b,c)$, then the word obtained from $v$ by deleting the letters other than $\mathcal{A}_c^{\pm 1} \cup \{\zeta^{\pm 1}\}$ also represents $\zeta^{r_i}$. And according to 2.1(3), if this redacted word has length $m_i$ then $r_i \leq Km_i^c$.

The $\zeta$-arcs in the boundary of the face of type $\mathcal{R}_c$ which we are considering span *disjoint* arcs in $\partial \Delta'$, so by applying the argument of the preceding paragraph to these disjoint arcs we get:

$$\sum_{i=1}^s |r_i| \leq \sum_{i=1}^s Km_i^c \leq Km^c.$$

□

One obtains a retraction of $\Gamma(a,b,c)$ onto $G_a \times G_b$ by killing the generators $\mathcal{A}_c^{\pm 1} \cup \{\zeta\}$. Using this retraction in place of the one used in the proof of (4.1), and noting that $a \geq b$, we obtain:

**4.2 Lemma.** *If the $z$-arcs in the boundary of a face in $\Delta'$ of type $\mathcal{R}_a$ or $\mathcal{R}_b$ are labelled $z^{r_1}, \ldots, z^{r_s}$, then $\sum_i |r_i| \leq K\mu^a$, where $\mu$ is the number of occurences of letters from $\mathcal{A}_a^{\pm 1} \cup \mathcal{A}_b^{\pm 1} \cup \{z^{\pm 1}\}$ in $w$ (the boundary label of $\Delta'$).*

**Step 2:** *Bounding the sum of the areas of faces of type $\mathcal{R}_a$ and $\mathcal{R}_c$.*

In this step of the proof we bound the number of 2-cells necessary to fill the faces of $\Delta'$ which are of type $\mathcal{R}_a$ or $\mathcal{R}_c$. In fact, a sufficiently sharp bound can be obtained simply by estimating the area of the words labelling their boundaries, thought of as null-homotopic words for $G_a$ or $G_c$.

**4.3 Lemma.** *The sum of the areas of the words bounding faces of type $\mathcal{R}_a$ in $\Delta'$ is bounded above by:*
$$2Kn^{a+1}.$$

*Similarly, the sum of the areas of the words bounding faces of type $\mathcal{R}_c$ is bounded above by:*
$$2Kn^{c+1}.$$

*Proof.* Consider a face $E_j$ of $\Delta'$ that corresponds to a monochromatic region of type $\mathcal{R}_a$ in $\Delta$. $\partial E_j$ consists of a number of $z$-arcs interspersed with subarcs of $\partial \Delta'$, so the word labelling $\partial E_j$ is of the form $W_{E_j} = z^{r_{j,1}} U_{j,1} \ldots z^{r_{j,s}} U_{j,s}$. At the cost of applying at most $(\sum_i |U_{j,i}|)(\sum_i |r_{j,i}|)$ relators, we can commute the subwords $U_{j,i}$ past the subwords $z^{r_i}$, and then freely reducing we get a null-homotopic word of the form $W'_{E_j} = z^{\rho_j} U_j$, with $|\rho_j| \le \sum_i r_{j,i}$ and $U_j = U_{j,1} \ldots U_{j,s}$ in the free group on $\mathcal{A}_a \cup \{z\}$. But now, by 2.1(2), with respect to the natural presentation for $G_a$, we have $Area(W'_{E_j}) \le K|U_j|^{a+1}$. Thus,

$$Area(W_{E_j}) \le \big(\sum_i |U_{j,i}|\big)\big(\sum_i |r_{j,i}|\big) + K\big(\sum_i |U_{j,i}|\big)^{a+1}.$$

We then use (a weak form of) Lemma 4.2 to replace the second factor in the first summand:

$$Area(W_{E_j}) \le \big(\sum_i |U_{j,i}|\big) K n^a + K\big(\sum_i |U_{j,i}|\big)^{a+1}.$$

(Recall that $n = |w| = |\partial \Delta'|$.) By summing over all $j$ and noting that $\sum_{i,j} |U_{i,j}| \le |w|$, we obtain the desired bound on the areas of faces of type $\mathcal{R}_a$. The argument for faces of type $\mathcal{R}_c$ is entirely similar. $\square$

**Step 3 (the hardest):** *Bounding the area of each face of type $\mathcal{R}_b$.*

Lemma 1.4 does not yield an easy estimate on the number of monochromatic regions of type $\mathcal{R}_b$ in $\Delta$, and this makes the job of estimating the sum of the areas of the faces of type $\mathcal{R}_b$ in $\Delta'$ more difficult than in the case of faces of type $\mathcal{R}_a$ and $\mathcal{R}_c$.

We decompose $\Delta'$ as the union of those subdiscs which are the closures of the connected components of $\Delta'$ minus the union of its $\zeta$-arcs. Some of these subdiscs are faces of type $\mathcal{R}_c$ — we are not concerned with these. We fix our attention on a subdisc $D_j$ which is not of this type. Notice that no edge in $\partial D_j \cap \partial \Delta'$ is labelled by a letter from $\mathcal{A}_c$, for otherwise the sub-diagram of $\Delta$ corresponding to $D_j$ would contain a monochromatic region of type $\mathcal{R}_c$ as a proper sub-diagram, and hence would have a $\zeta$-arc in its interior.

Let $\mu_j$ be the number of edges in $\partial D_j$ that are labelled by an element of $\mathcal{A}_a^{\pm 1} \cup \mathcal{A}_b^{\pm 1} \cup \{z^{\pm 1}\}$. As in (4.2), the sum of the lengths of the $z$-arcs in the boundary of each face of $D_j$ is at most $K\mu_j^a$. Consider a face in $D_j$ that is of type $\mathcal{R}_b$: its boundary consists of a number of $z$-arcs (total length $\le K\mu_j^a$), a number of $\zeta$-arcs (total length $\le Kn^c$) and a number of sub-arcs of $\partial \Delta'$ (total length $\le \mu_j$). We denote this face $E'$.

According to (2.1(3)), if a $z$-arc has length $r_l$, then one can replace it by a word of length at most $Kr_l^{1/b}$ in the generators $\mathcal{A}_b^{\pm 1} \cup \{z^{\pm 1}\}$, and by (2.1(2)) we know that this may be done at the cost of applying at most

$$(1) \qquad \alpha_b(Kr_l^{1/b})^{b+1} \le K^{2+(1/b)} r_l^{1+(1/b)}$$

relators. We shall replace each of $z$-arcs in $\partial E'$ in this way. We count how many relators we must apply in order to do so:

Let $\tilde{r}_l$ denote the number of occurences of letters from $\mathcal{A}_a^{\pm 1} \cup \mathcal{A}_b^{\pm 1} \cup \{z^{\pm 1}\}$ in the sub-arc of $\partial D_j$ corresponding to $z^{r_l}$. We have $r_l \leq K\tilde{r}_l^a$. And $\sum_l \tilde{r}_l \leq \mu_j$. These estimates, together with (1), show that at the cost of applying at most

$$\sum_l K^{2+(1/b)} r_l^{1+(1/b)} \leq K^{2+(1/b)} \left(\sum_l r_l\right)^{1+(1/b)}$$

$$\leq K^{2+(1/b)} (K\mu_j^a)^{1+(1/b)}$$

(2) $$=: K' \mu_j^{a+(a/b)}$$

relators, we can replace all of the labels on $z$-arcs in $\partial E'$ with equivalent words in the generators $\mathcal{A}_b^{\pm 1} \cup \{z^{\pm 1}\}$. The sum of the lengths of these equivalent words is at most

$$\sum_l K r_l^{1/b} \leq K \sum_l (K\tilde{r}_l^a)^{1/b}$$

$$= K^{1+(1/b)} \sum_l \tilde{r}_l^{(a/b)}$$

(3) $$\leq K' \mu_j^{(a/b)}.$$

At this stage, we have replaced the word labelling $\partial E'$ by a word which consists of at most $K'\mu_j^{(a/b)}$ letters from $\mathcal{A}_b \cup \{z\}$ (or their inverses) and some number of $\zeta$ (or its inverse), depending on how many of the $\zeta$-arcs in $\partial D_j$ are contained in $\partial E'$; let us write $Z(E')$ for the exact number of occurences of $\zeta$ or its inverse. (Later, when we sum over all $E'$ in $D_j$, we shall use the fact that $\sum Z(E') \leq Kn^c$.) Since the presentation of $\Gamma(a,b,c)$ with which we are working contains the relations $[\zeta, x] = 1$ for all $x \in \mathcal{A}_b \cup \{z\}$, we may commute all of the other letters past all occurences of $\zeta^{\pm 1}$ at a cost of applying at most

(3) $$K' \mu_j^{(a/b)} Z(E')$$

relators. Thus (after freely reducing) we obtain a word of the form $U\zeta^p$, where $U$ is a word of length at most $K'\mu_j^{(a/b)}$ in the generators $\mathcal{A}_b \cup \{z\}$ and their inverses. But $U\zeta^p$ represents the identity in $G_b \times \langle\zeta\rangle$, hence $p = 0$ and, by 2.1(2),

(4) $$Area(U) \leq K|U|^{b+1} \leq K(K')^{b+1} \mu_j^{a+(a/b)}.$$

Combining estimates (2) through (4), we see that for an appropriate constant $K''$ the following is an upper bound on the area of the word labelling $\partial E'$:

$$K'' \left(\mu_j^{a+(a/b)} + \mu_j^{(a/b)} Z(E')\right).$$

Thus, summing over all faces of type $\mathcal{R}_b$ in $D_j$, we obtain a *total* area of at most:

$$K'' \left(N_j \mu_j^{a+(a/b)} + \mu_j^{(a/b)} n^c\right),$$

where $N_j$ is the number of faces of type $\mathcal{R}_b$ in $D_j$.

Finally, summing over all $D_j$, and using the fact that $\sum \mu_j \leq n$ (because this is a sum of lengths of disjoint arcs in $\partial \Delta'$), we obtain:

**4.4 Lemma.** *The sum of the areas of the words bounding faces of type $\mathcal{R}_b$ in $\Delta'$ is at most*

$$K'' \left( N\ n^{a+(a/b)} + n^{c+(a/b)} \right),$$

*where $N$ is the number of faces of type $\mathcal{R}_b$.*

**Step 4:** *Bounding the number of faces of type $\mathcal{R}_b$ in $\Delta'$.*

By combining Lemma 4.3 and 4.4, and recalling that $1 \le b \le a < c$, we see that Theorem 2.2 is a consequence of the following assertion:

**4.5 Lemma.** *$\Delta'$ contains at most $n$ faces of type $\mathcal{R}_b$.*

*Proof.* Consider the finite connected graph $\mathcal{G}$ whose vertices are the faces of $\Delta'$, and which has an edge connecting two faces if and only if they abut along a $z$-arc or a $\zeta$-arc. We label the vertices $a, b$ or $c$, according to the type of face in $\Delta'$ to which they correspond. Note that $\mathcal{G}$ is a tree, because no $z$-arcs or $\zeta$-arcs intersect in the interior of $\Delta'$. Every vertex which is adjacent to a vertex of type $a$ or $c$ must be of type $b$, and vice versa.

We fix a vertex $x_0$ of valence one in $\mathcal{G}$, and consider the map which sends each vertex $x \in \mathrm{Vert}(\mathcal{G}) - \{x_0\}$ to the first vertex along the unique arc joining $x$ to $x_0$ in $\mathcal{G}$. This map sends the set of vertices that are of type $a$ or $c$ onto the set of vertices of type $b$ that have valence greater than one. But at most $n$ vertices of $\mathcal{G}$ correspond to faces of $\Delta'$ that contain an edge, and this collection of vertices includes all those of type $a$ and $c$ and all those of valence one. Thus $\mathcal{G}$ has at most $n$ vertices of type $b$, and hence $\Delta'$ has at most $n$ faces of type $\mathcal{R}_b$. □

## 5. Further Examples.

Thus far we have concentrated on the groups $\Gamma(a, b, c)$. We have done so in order to simplify the exposition as much as possible. However, the methods introduced here apply, more or less directly, to many related classes of groups. For example, a direct translation of the arguments presented above can be used to prove:

**5.1 Theorem.** *Let $a$ and $b$ be positive integers and let $G_a$ and $G_b$ be as defined in Section 2. Consider the group $J(a, b)$ obtained by amalgamating two copies of $G_a$ and $G_b$ along central cyclic subgroups in the following manner:*

$$J(a, b) := G_a *_{z_a = z_b} (G_b \times G'_b) *_{z'_a = z'_b} G'_a.$$

*If $a \ge b$ and $(a/b)^2 > a + 1$, then the Dehn function of $J(a, b)$ is $\simeq n^{(a/b)^2}$.*

Equally, instead of using the groups $G_a$ as basic building blocks, one could use groups of the form $\mathbb{Z}^a \rtimes F$, where $F$ is a non-abelian free group with a basis all of whose elements act by the unipotent matrix $\phi_a$ which has 1's on the diagonal and super-diagonal and zeros elsewhere.

The interested reader should be able to construct other examples by focusing on the key role played by subgroup distortion, cf. (3.4).

**Section 6: The proof of Theorem B.**

The following concludes our consideration of subgroup distortion.

**Theorem B.** *For every positive rational number $r$ there exist pairs of finitely presented groups $H \subseteq G$ with distortion $\simeq n^r$.*

*Proof.* We saw in (3.3) that if $a \geq b$ then the distortion $\delta_{a,b}(n)$ of $G_b$ in $G_a *_{\langle z \rangle} G_b$, the group formed by amalgamating $G_a$ and $G_b$ along their centres, is $\succeq n^{\frac{a}{b}-1}$. On the other hand, the proof of Proposition (3.1) shows that $K n^{c+1} \delta_{a,b}(n)$ is a lower bound for the Dehn function of $\Gamma(a,b,c)$, and we now know that this Dehn function is actually $\simeq n^{c+\frac{a}{b}}$ if $c > a$. Thus $\delta_{a,b}(n) \simeq n^{\frac{a}{b}-1}$. $\square$

## 7. Isodiametric Inequalities.

Isoperimetric inequalities measure the complexity of the word problem in a group by bounding the number of relators that one must apply in order to show that a null-homotopic word represents the identity. An alternative way of measuring complexity is to bound the length of the elements by which one must conjugate the relators being applied. Geometrically, this corresponds to bounding the diameter of van Kampen diagrams rather than their area. In certain contexts such bounds are the most appropriate measure of complexity (see [GS], [Ge2]).

**7.1 Definition.** Let $\langle \mathcal{A} \mid \mathcal{R} \rangle$ be a finite presentation. Let $w$ be a word that is null-homotopic with respect to this presentation and let $\Delta$ be a van Kampen diagram for $w$. Let $p$ be the basepoint of $\Delta$. Endow the 1-skeleton of $\Delta$ with a path metric $\rho$ that gives each edge length 1. The diameter of $w$ is defined by

$$diam(w) := \max\{\rho(p,q) \mid q \text{ a vertex of } \Delta\}.$$

The *isodiametric function* of $\langle \mathcal{A} \mid \mathcal{R} \rangle$ is

$$\Phi(n) := \max_{|w| \leq n} diam(w).$$

It is not hard to show that if two finite presentations define isomorphic (or quasi-isometric) groups then their isodiametric functions are $\simeq$ equivalent, and thus it makes sense to talk of the ($\simeq$ class of) the isodiametric function of a group.

*Terminology.* The term 'isodiametric function' is due to Gersten [Ge1]. Some authors prefer the term 'filling radius', which was introduced by Gromov [Gr2].

The purpose of this section is to prove Theorem C. In order to do so we consider the groups $J(a,b)$ defined in (5.1). By combining the usual presentations for $G_a$ and $G_b$ in the obvious way we obtain presentation for $J(a,b)$. More specifically, $G_a$ has generators

$$x_1, \ldots, x_{a-1}, z(=z_a), t_a$$

and defining relations

$[z, x_i] = [x_i, x_j] = 1 \ \forall i,j;\ [z, t_a] = 1;\ [x_i, t_a] = x_{i+1}$ if $i < a-1$ and $[x_{a-1}, t_a] = z$.

We denote this presentation $\langle \mathcal{A}_a, z \mid \mathcal{S}_a \rangle$. We write the corresponding presentations of $G'_a, G_b$ and $G'_b$ by ammending the subscript and adding primes appropriately. The natural presentation for $J(a,b)$ can then be written:

$J(a,b) = \langle \mathcal{A}_a, \mathcal{A}_b, \mathcal{A}'_a, \mathcal{A}'_b, z, z' \mid \mathcal{S}_a, \mathcal{S}_b, \mathcal{S}'_a, \mathcal{S}'_b;\ [b, b'] = 1 \ \forall b \in \mathcal{S}_b \cup \{z\}\ \forall b' \in \mathcal{S}'_b \cup \{z'\} \rangle.$

Henceforth we work with the above presentations of $G_a$ and $J(a,b)$.

**7.2 Theorem.** *If $a \geq b$ then the isodiametric function of $J(a,b)$ is $\simeq n^{a/b}$.*

Most of the ideas needed to prove this theorem are already present in the proof of Theorem 2.2, so we shall try to avoid repetition as much as possible.

In place of 2.1(2) we have:

**7.3 Lemma.** *There exists a constant $C_a > 1$ such that $diam(w'z^n) \leq C_a |w'|$ for all $n \in \mathbb{Z}$ and all words $w = w'z^n$ that represent the identity in $G_a$.*

*Proof.* A standard argument (see [B1] fig. 2 for example) shows that if a group admits a bounded asynchronous combing and $u$ is a combing line, then for any null-homotopic word $vu$ one has $diam(vu) \leq k|v|$, where $k$ is the fellowtraveller constant of the combing. $z^n$ is a combing line in the asynchronous combing of $G_a$ constructed in [B2]. □

**7.4 Corollary.** *Let $W$ be a word that represents the identity in $G_a$. Let $n_1$ be the number of occurences of $z^{\pm 1}$ in $W$ and let $n_2 = |W| - n_1$. There is a van Kampen diagram for $W$ in which every vertex is a distance at most $C_a n_2$ from the boundary.*

*Proof.* Given any group presentation with generators including $\mathcal{B} \cup \{\zeta\}$ and relations including $[b, \zeta] = 1 \,\forall b \in \mathcal{B}$, one can rewrite an arbitrary word $U$ in the form $U'\zeta^r$ simply by applying the relations $[b, \zeta] = 1$. This gives a van Kampen diagram with boundary label $U\zeta^{-r}U'^{-1}$ that can be drawn on the rectangular tessalation of the plane. Every vertex of this diagram can be connected to the arc of the boundary labelled $U$ by an edge-path (with all labels from $\mathcal{B}$) of length at most $|U'|$.

In order to obtain the desired diagram for $W$ one takes the diagram described above and attaches to it the van Kampen diagram for $W'z^{n_1}$ given by lemma (7.3). □

In place of (3.4) we have the following result, which is proved by an easy diagrammatic argument.

**7.5 Lemma.** *Let $\langle \mathcal{A} \mid \mathcal{R} \rangle$ and $\langle \mathcal{A}' \mid \mathcal{R}' \rangle$ be group presentations. Consider their cartesian product $\langle \mathcal{A}, \mathcal{A}' \mid \mathcal{R}, \mathcal{R}', [a, a'] = 1 \,\forall a \in \mathcal{A} \,\forall a' \in \mathcal{A}' \rangle$. For all words $U \in F(\mathcal{A})$ and $U' \in F(\mathcal{A}')$ one has:*

$$diam\left(UU'U^{-1}U'^{-1}\right) \geq \frac{1}{2} \min\left\{d_{\mathcal{A}}(1, U), d_{\mathcal{A}'}(1, U')\right\}.$$

We shall also need an inequality in the other direction (which is equally easy to prove):

**7.6 Lemma.** *Let $w$ be null-homotopic word for the product presentation in (7.5). Let $w_1$ be the word obtained by deleting from $w$ all occurences of the letters $\mathcal{A}'^{\pm 1}$ and let $w_2$ be the word obtained by deleting all occurences of the letters $\mathcal{A}^{\pm 1}$. Then,*

$$diam(w) \leq \max\left\{|w_1|, |w_2|, diam(w_1), diam(w_2)\right\}.$$

We shall obtain the lower bound needed for Theorem 7.2 by considering the sequence of words

$$W_n = V_n V_n' V_n^{-1} V_n'^{-1},$$

where $V_n$ is the shortest word in the generators $\mathcal{A}_a^{\pm 1}$ representing $z^{n^a}$ and $V_n'$ is the same word over the alphabet $\mathcal{A}_a'^{\pm 1}$.

**7.7 Lemma.** *With respect to the natural presentation of $J(a,b)$,*

$$diam(W_n) \geq \frac{1}{2} k_b n^{a/b},$$

*where $k_b$ is the constant of (2.1).*

*Proof.* We consider a minimal diameter van Kampen diagram for $W_n$. By arguing as in the proof of (3.1) one sees that this diagram must contain a subdiagram whose boundary is labelled $z^{n^a} z'^{n^a} z^{-n^a} z'^{-n^a}$.

One would like to say that this subdiagram is actually a diagram over the sub-presentation $\langle \mathcal{A}_b, \mathcal{A}'_b, z, z' \mid \mathcal{S}_b, \mathcal{S}'_b, [b,b'] = 1 \, \forall b \in \mathcal{S}_b \cup \{z\}, \forall b' \in \mathcal{S}'_b \cup \{z'\}\rangle$. If one could legimately do so, then an appeal to Lemmas 7.3 and 2.1(3) would complete the proof.

However, minimal diameter van Kampen diagrams can contain subdiagrams in their interior that are bounded by loops whose labelling word freely reduces to the empty word, and because of this it might be the case that the subdiagram under consideration contains monochromatic regions of type $\mathcal{S}_a$ or $\mathcal{S}'_a$. (The boundary cycle of such a region would be labelled by a word in $z$ or $z'$ that freely reduced to the empty word.) This problem is easily remedied: we simply excise all such regions and replace them with diagrams of zero area. (This does not increase the diameter of the ambient diagram.) □

It remains to establish the upper bounds needed for Theorem 7.2. Let the constants $C_a$ and $C_b$ be as in (7.3) and let $K_a$ and $K_b$ be as in (2.1(3)). Let $\kappa = 2C_a C_b K_a K_b$.

**7.8 Proposition.** *In the natural presentation of $J(a,b)$, every null-homotopic word satisfies $diam(w) \leq \kappa |w|^{a/b} + |w|$.*

*Proof.* Let $n = |w|$. It suffices to consider the case where $w$ is irreducible. As in Section 4, we begin with a van Kampen diagram $\Delta$ for $w$ and abstract from it only the pattern of $z$-arcs and $z'$-arcs (as in fig.1). Continuing the analogy with Section 4, we shall assemble a new diagram for $w$ by filling the faces of this simplified version of $\Delta$ with disc subdiagrams: each subdiagram will be a van Kampen diagram over one of the natural sub-presentations of $J(a,b)$ corresponding to the subgroups $G_a, G'_a$ or $G_b \times G'_b$. We shall refer to these different types of subdiagram as type $a, a'$ and $b$ respectively.

The subdiagrams are attached along $z$-arcs and $z'$-arcs. The boundary of each subdiagram contains at most $K_a n^a$ edges, at most $n$ of which come from $\partial \Delta$. Subdiagrams of type $a$ and $a'$ can only abut subdiagrams of type $b$.

We number the remaining paragraphs of the proof for ease of reference.

*7.8.1:* It suffices to fill the subdiagrams in such a way that every vertex is a distance at most $\kappa n^{a/b}$ from the boundary of the ambient diagram.

*7.8.2:* We fill the faces of type $a$ and $a'$ with the diagrams furnished by (7.4). Each vertex in the resulting subdiagrams can be connected to the boundary of its subdiagram by an edge path of length at most $K_a n$.

*7.8.3:* We now begin filling in the faces of type $b$, beginning with a neighbourhood of the $z$-arcs in their boundaries. Given a $z$-arc labelled $z^r$, we choose a geodesic word $u$ representing $z^r$ in $G_b$ and attach to the $z$-arc a van Kampen diagram for $uz^{-r}$ over $\langle \mathcal{A}_b, z \mid \mathcal{S}_b \rangle$. This diagram is chosen as in Lemma 7.4.

*7.8.4:* With this choice, we can connect any vertex on the $z$-arc under consideration to the initial point of the arc labelled $u$ by a path of length at most $C_a|u|$. We must estimate $|u|$. ¿From 2.1(3) we have $|u| \le K_b r^{1/b}$, and we also know $r \le K_a n^a$. Thus $|u| \le K_a K_b n^{a/b}$. (This last bound is rather crude. As in (4.2), we can replace the estimate $r \le K_a n^a$ by $\sum r \le K_a n^a$, where the sum is taken over all $z$-arcs in the boundary of a fixed subdiagram of type $b$.)

*7.8.5:* We apply the construction of the previous paragraph to every $z$-arc and every $z'$-arc in the boundary of each face of type $b$. By doing so, we replace the original boundary label of the face by a word of the form $u_1 u'_1 \ldots u_m u'_m$, where the $u_i$ are words in the generators $\mathcal{A}_b \cup \{z\}$ and the $u'_i$ are words in the generators $\mathcal{A}'_b \cup \{z'\}$. The paranthetical remark at the end of the last paragraph tells us that both $\sum |u_i|$ and $\sum |u'_i|$ are bounded above by $K_a n^a$. And Lemma 7.3 (for $G_b$ and $G'_b$) then implies that $diam\,(u_1 \ldots u_m)$ and $diam\,(u'_1 \ldots u'_m)$ are bounded by $C_b K_a n^a$. We may therefore fill the remaining part of each face of type $b$ with the diagram of diameter at most $C_b K_a n^a$ provided by Lemma 7.6 (taking the basepoint of the diagram on $\partial \Delta$).

*Claim.* Every vertex of the van Kampen diagram just constructed is within a distance $\kappa n^{a/b}$ of the boundary and hence the diagram has diameter at most $\kappa n^{a/b}+n$).

*7.8.6:* First consider the case of a vertex in a subdiagram of type $a$ or $a'$: it is within a distance $K_a n$ of some point $p$ on the boundary of the subdiagram (7.8.2). If $p$ lies on $\partial \Delta$ then we are done. If not then $p$ lies on a $z$-arc, and on the other side of this $z$-arc we have attached a van Kampen diagram over $G_b$ or $G'_b$ (see 7.8.3). The attached diagram has its basepoint on $\partial \Delta$ and it has diameter $\le C_a K_a K_b n^{a/b}$ (see 7.8.4). Thus we may connect our original vertex to $\partial \Delta$ by a path of length at most $K_a n + C_a K_a K_b n^{a/b}$, which is less than $\kappa n^{a/b}$.

*7.8.7:* Finally, consider the case of a vertex in a subdiagram of type $b$. If this vertex lies in one of the subdiagrams that we attached to $z$-arcs in (7.8.3), then it is a distance $\le C_a K_a K_b n^{a/b}$ from the boundary of the whole diagram (7.8.4). On the other hand, if it lies in the remainder of the subdiagram then by (7.8.5) it is a distance at most $C_b K_a n^a$ from the boundary. $\square$

*Acknowledgement.* I am grateful to Dani Wise and Hamish Short, whose detailed comments improved the exposition of this article.

MATHEMATICAL INSTITUTE, 24–29 ST. GILES, OXFORD OX1 3LB, GREAT BRITAIN
*E-mail address*: bridson@math.princeton.edu